\documentclass[11pt]{article}
\usepackage{bbm}
 \usepackage{amssymb}
\usepackage{amssymb, amsthm, amsmath, amscd}
\usepackage[usenames,dvipsnames]{color}
\usepackage{hyperref}
\hypersetup{hidelinks, colorlinks=true, allcolors=blue, pdfstartview=Fit, breaklinks=true}

\newtheorem{thm}{Theorem}[section]

\newtheorem{cor}[thm]{Corollary}
\theoremstyle{definition}
\newtheorem{NN}[thm]{}
\theoremstyle{definition}
\newtheorem{df}[thm]{Definition}
\theoremstyle{definition}
\newtheorem{rem}[thm]{Remark}

\theoremstyle{definition}

\newcommand{\red}{\textcolor{red}}
\newcommand{\blue}{\color{blue}}

\renewcommand{\phi}{\varphi}

\newcommand{\aff}{\rm aff}

\newcommand{\N}{\mathbb{N}}
\newcommand{\Z}{\mathbb{Z}}
\newcommand{\Q}{\mathbb{Q}}
\newcommand{\R}{\mathbb{R}}
\newcommand{\C}{\mathbb{C}}
\newcommand{\T}{\mathbb{T}}

\numberwithin{equation}{section}

\newcommand{\Aff}{\operatorname{Aff}}

\newcommand{\id}{\operatorname{id}}

\newcommand{\hm}{homomorphism}
\newcommand{\dt}{\delta}
\newcommand{\ep}{\varepsilon}

\newcommand{\td}{\tilde}

\newcommand{\lr}{\longrightarrow}
\newcommand{\ld}{\lambda}
\newcommand{\cd}{\cdots}

\newcommand{\DT}{\Delta}

\newcommand{\la}{\langle}
\newcommand{\ra}{\rangle}
\newcommand{\andeqn}{\,\,\,{\rm and}\,\,\,}
\newcommand{\rforal}{\,\,\,{\rm for\,\,\,all}\,\,\,}
\newcommand{\CA}{$C^*$-algebra}
\newcommand{\SCA}{$C^*$-subalgebra}

\newcommand{\af}{{\alpha}}
\newcommand{\bt}{{\beta}}

\definecolor{purple}{RGB}{150,10,200} 
\newcommand{\pcl}{\color{purple}}

\newcommand{\beq}{\begin{eqnarray}}
\newcommand{\eneq}{\end{eqnarray}}
\newcommand{\tforal}{\,\,\,\text{for\,\,\,all}\,\,\,}

\newcommand{\zo}{{\cal Z}_0}

\newcommand{\LAff}{{\rm LAff}}

\usepackage{amsfonts}
\usepackage{mathrsfs}
\usepackage{textcomp}
\usepackage[all]{xy}

\title{A review of the Elliott program of classification
of simple amenable \CA s}

\author{Guihua Gong, Huaxin Lin,  and Zhuang Niu
 }
\date{
}

\begin{document}

\maketitle

\begin{abstract}
We give a brief survey of the development of the Elliott program of
classification of separable simple amenable \CA s.  
\end{abstract}

\section{Introduction}

By a \CA, we mean a norm closed and adjoint closed subalgebra
of the algebra $B(H)$ of all bounded linear operators 
on some Hilbert space $H.$
There is an abstract definition of Gelfand and Naimark (\cite{GN}) for \CA s, namely,
a Banach algebra with an involutive $*$-operation such that $\|x^*x\|=\|x\|^2$
for all $x\in A.$
Let $X$ be a compact metric space, $G$ be a group, $\af$ be a \hm\,
from $G$ to the group of homeomorphisms of $X$ and $\mu$ be a faithful
$\af$-invariant
Borel probability measure. We have  a representation
$\pi: G\to U(B(H)),$ the unitary group of $B(H),$ where $H=L^2(X, \mu).$
For each $f\in C(X),$ we have a bounded operator $M_f$ defined by $M_f(g)=fg$
for all $g\in L^2(X, \mu).$  Then we obtain a \CA\, $C(X)\rtimes_\af G,$
the \CA\, generated by $\pi(G)$ and $\{M_f: f\in C(X)\}$ as a \SCA\, of $B(H).$ This \CA\, is often used to
study the dynamical system $(X, G, \af).$
\CA s may also arise from non-commutative geometry, number theory and group representations,
to name a few related fields.   Needless to say, it is of paramount importance to determine the structure of a given \CA\, from a few known  and possibly computable data and
to recognize  that some \CA s from different subjects may  be isomorphic as \CA s.
This leads  to the question of  classification of  \CA s.

In the First International Congress of Basic Science,  David Mumford
gave the opening  plenary lecture on Consciousness, robots and DNA (\cite{Mf}).
At the very beginning, he cautioned the audience  that there would be no mathematics in his talk.
Nevertheless, by the very end of this magnificent lecture, when he ended his talk
on DNA being a measuring instrument opens a Pandora's box, he proposed a problem in
\CA s about non-commuting self-adjoint operators to be approximated by commuting ones,
which is closely related to the Elliott program of classification of amenable \CA s (see \cite{Linmatrices}).

Classification of amenable \CA s might begin with Glimm's classification of UHF-algebras (by supernatural numbers) in the late 1950's. This was followed by Dixmier's
      classification of matroid \CA s in the 1960's. Ola Bratteli used the Bratteli diagram to classify AF-algebras,
      i.e., \CA s which are inductive limits of finite dimensional \CA s.  George A. Elliott then (1976)
      used dimension groups to classify AF-algebras (\cite{Ellaf}).
      It was later realized   that
      the dimension groups are in fact the $K_0$-groups. Even though $K_1$-groups of AF-algebras
      are trivial, realizing dimension groups are $K_0$-groups is very important since, immediately,
      $K_1$-groups enter the picture.

      By early 1990, Elliott proposed
      to classify simple \CA s by their K-theory.

    As suggested by    Elliott (\cite{Ellicm}),
separable simple
\CA s $A$ are
naturally divided into three  cases
according to their
 $K_0$-groups:\\
{{Case 1: $K_0(A)\not=\{0\}$ is a (directed) ordered group,}}\\
{{Case 2.
$K_0(A)_+=\{0\}=V(A)$,}}\\
Case 3: $K_0(A)=K_0(A)_+$, and $V(A)\not=\{0\}$,\\
where $V(A)$ is the set of Murry-von Neumann  equivalence classes of projections in {{the stabilization $A\otimes {\cal K}.$

For  Case 3, i.e., the case of purely infinite simple \CA s, the classification by the Elliott invariant
was accomplished in the last century by the work of Kirchberg and Phillips (see \cite{KP} and \cite{Pclass})
with some earlier works of the M. R\o rdam  and others (cf.~\cite{R1}, \cite{R2},  \cite{LP1}, \cite{ERr} and
\cite{R3}).
This report will only discuss the finite case.

To clarify what classification of a class of \CA s  by the Elliott invariant is 
about, let us end this section  with the following statement.

For a complete classification of a class of \CA s by the Elliott invariant, there must be three  ingredients:

(1) A definition of invariants for \CA s in the class;  (2)  A (range) theorem  that says exactly what values of 
the invariants defined in (1) occur for \CA s in the class; (3) the  isomorphism theorem, i.e., 
the statement that any two \CA s in the class with  the same, i.e., isomorphic  invariants
are isomorphic.

%


      \section{The beginning}

      At the beginning of the 1990's,  Elliott  presented the following
      classification results:

      \begin{thm}[G. A.  Elliott \cite{Ellcl1}] \label{TElliott1}
      Let $A=\lim_{n\to\infty} A_n$ and $B=\lim_{n\to\infty}B_n$ be unital simple \CA s  with real rank zero (see \ref{Drealstable} below),
      where $A_n$ and $B_n$ are finite direct sums of \CA s of the form
      $M_{r(n,i)}(C(\T)),$ i.e., matrices over the continuous functions on the unit circle.

      Then $A\cong B$ if and only if
      $$
      (K_0(A), K_0(A)_+, [1_A], K_1(A))\cong (K_0(B), K_0(B)_+, [1_B], K_1(B)).
      $$
     \end{thm}

     \begin{thm}[G. A.  Elliott \cite{Ellcl1}]
      Let $G_0$ be a simple countable (ordered) dimension  group, not equal to $\Z,$
      with an order unit $u,$
      and let $G_1$ be a countable abelian group.
        It follows that there is a unital simple \CA\,
        %
       of real rank zero which is the inductive limit of a sequence of finite direct sums of
       circle algebras such that
       $$
       (K_0(A), K_0(A)_+, [1_A], K_1(A))=(G_0, (G_0)_+, u, G_1).
       $$
     %
     \end{thm}

      This is the beginning of an age that classes of simple \CA s became  classifiable.

      \begin{df}\label{Drealstable}
      Let $A$ be a \CA, $S\subset A$ be a subset, and $x\in A$ be an element.
      Let $\ep>0.$
      We write $x\in_\ep S,$ if there is $s\in S$ such that $\|x-s\|<\ep.$

      Denote by $T(A)$ the space of all tracial states of $A.$
      Denote by $\Aff(T(A))$ the space of all real affine continuous functions on $T(A).$

      A unital \CA\, $A$ is said to have stable rank one,
      if the set of invertible elements of $A$ is dense in $A.$
      A non-unital \CA\, has stable rank one if $\tilde A,$ the minimum
      unitization of $A,$ has stable rank one.

      A unital \CA\, $A$ is said to have real  rank zero,
      if the set of invertible self-adjoint elements of $A$ is dense in $A_{s.a.},$
      the self-adjoint part of $A.$
      A non-unital \CA\, has real  rank zero,  if $\tilde A$
      has real  rank zero.
      \end{df}

      A unital commutative \CA\, $A=C(X),$ where $X$ is a compact Hausdorff space,
      has stable rank one if and only if  ${\rm dim}X\le 1,$ and $A$ has real rank zero,
      if and only if ${\rm dim}X=0.$  Every von Neumann algebra has real rank zero.
      Every AF-algebra has real rank zero.   It was proved by S. Zhang \cite{Zh} that
      every purely infinite simple \CA\, has real rank zero.  It is then natural to begin
      the classification program  with
      simple \CA s of real rank zero  in mind such as  Theorem \ref{TElliott1}.

      \begin{df}\label{DAH}
      Denote by ${\cal I}_0$ the class of all finite dimensional \CA s.
      For $k\in \N,$
      denote by ${\cal I}_k$ the class of \CA s which are finite direct sums of \CA s
      with the form
      $
      PM_r(C(X))P,
      $
      where $X$ is a finite CW complex with
      ${\rm dim}(X)\le k$   and $P\in M_r(C(X))$ is a projection.

      A \CA\, $A$ is said to be an AH-algebra if
      $A=\lim_{n\to\infty} A_n,$ where each $A_n\in {\cal I}_{k(n)}$ for some $k(n)\in \N.$

      If there is $k\in \N$ such that $k(n)\le k, $ we say that the AH-algebra $A$ has no dimension growth.

      A unital simple \CA\ is said to be of
      {\em slow dimension growth}, if
      $A=\lim_{n\to\infty} A_n,$
      where $A_n=\bigoplus_{i=1}^{k(n)} P_{n,j}M_{r(n,j)}(C(X_{n,j})P_{n,j}$ and
      $P_{n,j}M_{r(n,j)}(C(X_{n,j})P_{n,j}\in \bigcup_k {\cal I}_k$ such that
      \beq
      \limsup_{n\to\infty} \max_{1\le j\le k(n)}{{\rm dim} (X_{n,j})\over{{\rm rank}(P_{n,j})+1}}=0.
      \eneq

         Let $I_k=\{f\in M_k(C[0,1]):~ f(0)=\ld_0\cdot 1_k,f(1)=\ld_1\cdot 1_k, \ld_0,\ld_1\in \C\}$.
         A \CA\, with the form  $M_n(I_k)$ is called a dimension drop $C^*$-algebras. 

      A \CA\, is said to be sub-homogeneous, if it is a \SCA\, of a \CA\, in ${\cal I}_k$ for some
      $k\in \N.$  A \CA\, is said to be an ASH-algebra if it is an inductive limit of sub-homogeneous \CA s.

      \end{df}

      \begin{df}
          Let $A$ be a unital separable simple \CA\, with tracial state space $T(A).$
      For each $\tau\in T(A),$ we extend it to a densely defined trace $\tau\otimes Tr$
      on $A\otimes {\cal K},$
      where $Tr$ is the densely defined lower semi-continuous trace on ${\cal K},$
      and ${\cal K}$ is the \CA\, of compact operators on $l^2.$ We continue to use
      $\tau$ for $\tau\otimes Tr.$
      Denote by $\rho_A: K_0(A)
      \to \Aff(T(A))$ the pairing defined by
      $\rho_A([p])=\tau(p)$ for all projections $p\in A\otimes {\cal K}$  and $\tau\in T(A).$

      Let $A$ be a unital simple exact \CA\,  and $p, q\in A\otimes {\cal K}$ be projections.
      We say  that $A$ has
      strict comparison for projections  if  $\tau(p)>\tau(q)$ for all $\tau\in T(A)$ implies
      that $q\lesssim p.$\footnote{We avoid the discussion of quasitraces here.}
      \end{df}

       The following theorem is a combination of proofs in 
       \cite{DNNP} and \cite{BDR} (see also  (iii) of Theorem 7.3.9 of  \cite{Bproj}).

      \begin{thm}[cf. D{\u a}d{\u a}rlat-Nagy-N\'{e}methi-Pasnicu \cite{DNNP} and
      Blackadar-D{\u a}d{\u a}rlat-R\o rdam \cite{BDR}]\label{DNNP}
      A unital simple AH-algebra with slow dimension growth has stable rank one
      and has 
      weakly unperforated $K_0(A).$
      \end{thm}

      After Elliott's theorem \ref{TElliott1} of classification of
      unital simple
      AH algebras with  real rank zero whose building blocks have spectrum $S^1$, Hongbing Su generalized Elliott's results to certain ASH 
      inductive limits with sub-homogeneous building blocks whose spectra are  locally one dimensional   but non-Hausdorff spaces
      (\cite{Su}). However,  to deal with AH inductive limit systems with higher dimensional spectrum,  one faced the difficulty of dealing with
       building blocks which are not semi-projective.
     A method
    to break the spectrum and  reduce the dimension (see \cite{Linnormal94}, \cite{Linnormal96} and \cite{EGLP}) was developed.
    Such a method involves  decomposing elements into two parts with the major parts being approximated by elements in a simpler $C^*$-algebra  such as a finite dimensional $C^*$-algebra, and plays an important role in the  classification theory which eventually led to the notion of tracial rank and generalized tracial rank.

      This approach led to an  important early breakthrough  in  the Elliott program.

      \begin{thm}[Elliott-Gong, Dadarlat-Gong \cite{EG1}, \cite{D}, \cite{G3-4} and  \cite{DG}]\label{TEG}
      Let $A$ and $B$ be unital simple AH-algebras with no dimension growth of real rank zero.
      Then $A\cong B$ if and only if
      \beq
       (K_0(A), K_0(A)_+, [1_A], K_1(A))\cong (K_0(B), K_0(B)_+, [1_B], K_1(B)).
      \eneq

      \end{thm}

      \begin{df}\label{DElliottinv}
      By an ordered group, we mean an abelian group which is (partially) ordered.
      Let $G$ be an ordered group. An element $u\in G$ is said to be an order unit
      if, for any $g\in G,$ there is $n\in \N$ such that  $g\le n u.$
      An ordered group is said to be simple if every $g\in (G)_+\setminus \{0\}$
      is an order unit.  A simple ordered group $G$ is said to be weakly unperforated if, for any $x\in G,$
$mx\in G^+\setminus \{0\}$ for some integer $m\in \N$ implies $x\in (G)_+\setminus \{0\}.$
      Denote by $S(G,u)$ the state space of $(G, u),$ i.e.,
      the set of all order preserving \hm s $s: G\to \R$ such that $s(u)=1.$
      It is a convex set.

     Denote by $G_{tor}$ the torsion subgroup of $G.$

       For any convex set $S,$ denote by $\partial_e(S)$ the set of extreme points of
       $S.$
      \end{df}

      \begin{df}
      Let $A$ be a unital separable simple \CA\, with tracial state space  $T(A).$
     %
      %
      Recall that $S(K_0(A), [1_A])$  is the state space of $K_0(A).$
    The map $\rho_A: K_0(A)\to \Aff(T(A))$ induces an affine map
     $\rho_A^\sharp:  T(A)\to  S(K_0(A), [1_A])$ by $\rho_A^\sharp(\tau)(x)=\rho_A(c)(\tau)$
     for all $\tau\in T(A)$ and $x\in K_0(A).$

      Let us denote by ${\rm Ell}(A)$ the sextuple
      $$
      (K_0(A), K_0(A)_+, [1_A], K_1(A), T(A), \rho_A).
      $$
     Suppose that $B$ is another unital separable simple \CA. We write
     $$
     (K_0(A), K_0(A)_+, [1_A], K_1(A), T(A), \rho_A)\cong(K_0(B), K_0(B)_+, [1_B], K_1(B), T(B), \rho_B),
     $$
     if there is an ordered group isomorphism
     $\kappa_0: K_0(A)\to K_0(B)$ with $\kappa_0([1_A])=[1_B],$ a group isomorphism
     $\kappa_1: K_1(A)\to K_1(B)$ and an affine homeomorphism
     $\kappa_T: T(A)\to T(B)$
     such that $\rho_B(\kappa_0(x))(\tau)=\rho_A(x)(\kappa_T^{-1}(\tau))$
     for all $x\in K_0(A)$ and $\tau\in T(B).$
     If we agree that $K_0(A)_+=\{x\in K_0(A): \rho_A(x)>0\}\cup \{0\},$
     then, we may write
     $$
     {\rm Ell}(A)=((K_0(A), [1_A], T(A), \rho_A), K_1(A)).
     $$
           \end{df}

       Based on the  \cite{Ell2}, \cite{Ell3}, \cite{NT}, \cite{Li1}, \cite{Li2}, \cite{Li3}, \cite{Li4},
      a further breakthrough was made by the following theorem:


      \begin{thm}[Elliott-Gong-Li \cite{G5} and \cite{EGL-AH}]\label{Tegl}
      Let $A$ and $B$ be unital simple \CA s.
      Suppose that both $A$ and $B$ are unital AH-algebras with no dimension
      growth. Then $A\cong B$ if and only if
      ${\rm Ell}(A)\cong {\rm Ell}(B).$
      \end{thm}

      Theorem \ref{Tegl} is made more important by the following range theorem:

      \begin{thm}[Villadsen \cite{V1}]\label{Trange}
      %
      Suppose that $G$  is a countable
      simple ordered group which is  weakly unperforated with Riesz interpolation such  that $G/G_{tor}$ is noncyclic,
      $u\in G^+,$  $H$ is a countable abelian group,  $\Delta$
       is a metrizable Choquet simplex, and
       $\rho: G\to \Aff(\Delta)$ is an order preserving \hm\, such that $\rho^\sharp(\partial_e(\Delta))=\partial_e(S(G, u)).$

        Then there exists a unital simple $AH$-algebra $A=\lim_{n\to\infty} A_n,$
        where each $A_n\in {\cal I}_3,$ such that
${\rm Ell}(A)=((G, u, H ),\Delta, \rho).$
      \end{thm}

       It should be noted that $\Z\oplus G_0,$ where $G_0$ is a torsion group,
       with $(\Z\oplus G_0)_+=\{(n,g): n\in \N,\, g\in G_0\}\cup \{(0,0)\}$
       is weakly unperforated but (for two reasons) excluded from the range described in Theorem \ref{Trange}.

      Fortunately for the further development of classification theory,  Jiang and Su provided the following Jiang-Su algebra ${\cal Z}.$ In particular,
      ${\cal Z}$ is not a unital simple AH-algebra.

      \begin{thm}[Jiang-Su \cite{JS}]\label{JS}
      There is an infinite dimensional simple \CA\, ${\cal Z}$ which is an inductive limit
      of dimension drop algebras such that
      $$
      {\rm Ell}({\cal Z})={\rm Ell}(\C).
      $$
      \end{thm}

      \begin{thm}[Gong-Jiang-Su \cite{GJS}]\label{GJS}
         For any unital simple $C^*$-algebra $A$, $(K_0(A\otimes {\cal Z}), K_0(A\otimes {\cal Z})_+)$ is weakly unperforated.  Furthermore,
      \beq
      {\rm Ell}(A\otimes {\cal Z})\cong {\rm Ell}(A)
      \eneq
      if and only if $(K_0(A), K_0(A)_+)$ is weakly unperforated. In particular, there are  unital simple AH algebras $A$  (such as the Villadsen algebra \cite{Vi}) which are not ${\cal Z}$-stable. 
      \end{thm}


      Next let us now introduce the \CA s of UCT class.

      \begin{df}[Rosenberg and Schochet  \cite{RS}]\label{DUCT}
      A  separable \CA\, $A$ is said to be in the UCT class, if it satisfies the following
      Universal Coefficient Theorem:
      For any $\sigma$-unital \CA\, $B,$
      one has the following short exact sequence:
      \beq
      0\to {\rm ext}_\Z^1(K_*(A), K_*(B))\, {\stackrel{\dt}{\to}}\,  KK^{*}(A, B)\, {\stackrel{\gamma}{\to}} \,{\rm Hom}(K_*(A), K_*(B))\to 0,
      \eneq
      where $\gamma$ has degree $0$ and $\dt$ has degree $1.$ The sequence
      is natural in each variable and splits unnaturally.

      Let ${\cal N}$ be the smallest class of separable amenable \CA s such that
      (1) $\C\in {\cal N};$ (2) ${\cal N}$ is closed under countable inductive limits,
      (3) If one has the short exact sequence
      $0\to A\to E\to B\to 0,$ and two of the terms are in ${\cal N},$ then so is the third;
      (4) ${\cal N}$ is closed under $KK$-equivalence (see Definition 19.1.1 of \cite{Bb} for the definition).

      It is known that ${\cal N}$ 
      is closed under stable isomorphism (by (4) above) and 
      contains all inductive limits of type $I$ \CA s (hence all ASH-algebras).
      It is closed under tensor products and crossed products by $\R$ and $\Z.$
      In particular, $O_n, O_\infty\in {\cal N}.$
      \end{df}

      \begin{thm}[Rosenberg and Schochet  \cite{RS}] \label{TUCT}
      Every separable \CA\,  in ${\cal N}$ is in the UCT class.
      \end{thm}

      Since the Elliott invariant includes $K$-theory, one expects that classes of separable
      \CA s to be classified by the Elliott invariant are  friendly to $KK$-theory. So we must also
      assume that the \CA s to be classified by the Elliott invariant are in the UCT class.

      \section{Abstract approach}

     Finite \CA s discussed in the previous section are presented as inductive limits of
     some basic building blocks such as finite dimensional \CA s, or some (sub)-homogeneous
     \CA s.  However, some naturally arising  simple \CA s, such as crossed product C*-algebras
      $C(X)\rtimes_\af \Z$ for minimal homeomorphisms $\af: X\to X$, do not have the appearance
     of inductive limits of some simpler  forms.

\begin{df}\label{Dher}
Let $A$ be a \CA.
Let  $a\in A_+.$  Denote by  ${\rm Her}(a)$ the hereditary \SCA\, $\overline{aAa}.$
If $a, b\in A_+,$ we write $a\lesssim b$ ($a$ is Cuntz subequivalent to $b$),
if there exists a sequence of elements $x_n\in A$ such that $a=\lim_{n\to\infty} x_n^*x_n$
and $x_nx_n^*\in {\rm Her}(b).$  {{If both $a\lesssim b$ and $b\lesssim a$
{{hold,}} then we {{write $a\sim b$ and}}
say $a$ is Cuntz equivalent to $b$.}} The Cuntz equivalence class represented by $a$ will
be denoted by $\la a\ra.$

\end{df}

  The following definition was inspired by the papers
  \cite{Gong-Lin-Math-scand}, \cite{EGLP}, \cite{Linnormal94},
     \cite{Linnormal96}, \cite{EG1}, \cite{G5}, and \cite{Ps}.

      \begin{df}[Lin, Definition 3.6.2 of \cite{Linbook}, also \cite{Lintrk} and \cite{LinTAF}]\label{Dtrk}
      Let $A$ be a unital simple  \CA.
      $A$ is said to have tracial rank at most $k,$ if, for any $\ep>0,$ any
      finite subset ${\cal F}\subset A$ and any given $a\in A_+\setminus \{0\},$
      there exists a \SCA\, $B\in {\cal I}_k$ with $p=1_B$ such that

      (i) $\|px-xp\|<\ep\rforal x\in {\cal F},$

      (ii) $pxp\in_\ep B\rforal  x\in {\cal F},$ and

      (iii) $p\lesssim a.$

      If $A$ has tracial rank at most $k,$ we will write $TR(A)\le k.$
      If $TR(A)\le k$ but $TR(A)\not\le k-1,$ we say that
      $A$ has tracial rank $k$ and write $TR(A)=k.$
      In particular, if $k=0,$ we say that $A$ has tracial rank zero and write
      $TR(A)=0.$

      If $TR(A)\le k$ for some $k\in \N,$ then $A$ has stable rank one and
      $K_0(A)$ is weakly unperforated.
       If $TR(A)=0,$ then $A$ also has real rank zero
        (see \cite{Lintrk} and \cite{LinTAI}).

       It was proved in \cite{EG1} that a unital simple  AH-algebra with   slow dimension growth
      and with real rank zero has tracial rank zero. 
      \end{df}

      \begin{thm}[Lin \cite{Linduke}]\label{Ltr=0}
      Let $A$ and $B$ be unital simple \CA s of tracial rank zero in the UCT class.
      Then $A\cong B$ if and only if
      \beq\nonumber
       (K_0(A), K_0(A)_+, [1_A], K_1(A))\cong (K_0(B), K_0(B)_+, [1_B], K_1(B)).
       \eneq
      \end{thm}

      It is proved in \cite{Linred} that
      a unital simple AH-algebra $A$ which has real rank zero, stable rank one, and weakly unperforated $K_0(A)$
      must have tracial rank zero.   In particular, together with Theorem \ref{DNNP}, all
      unital simple AH-algebras with slow dimension growth and real rank zero have tracial rank zero 
      (and  are
      therefore classifiable by Theorem \ref{Ltr=0} (compare Theorem \ref{TEG})). 
      It should be mentioned that 
      a unital separable simple exact \CA\, $A$ with real rank zero, stable rank one, and weakly unperforated 
      $K_0(A)$ has the property that $\rho_A(K_0(A))$ is dense in $\Aff(T(A))$ (see Lemma III 3.4 of 
      \cite{BH}).
      Many other simple \CA s have been proved to have tracial rank zero,
      for example,  in \cite{Wrr0}, W. Winter showed that a unital separable simple \CA\,  
      with finite decomposition rank and real rank zero has tracial rank zero (see \cite{NW} for further examples).
      For crossed product \CA s, let us mention the next result.

      \begin{thm}[Lin-Phillips, \cite{LP2}]\label{TLP2}
      Let $X$ be an infinite, compact finite dimensional metrizable space and
      $\af: X\to X$ be a minimal homeomorphism.
      Then the crossed product $A=C(X)\rtimes_\af\Z$
      has tracial rank zero if and only if
      $\rho_A(K_0(A))$ is dense in $\Aff(T(A)).$
            \end{thm}

      In particular, suppose that $X$ is a connected compact metric space with finite dimension,
      and $\af: X\to X$ is minimal and uniquely ergodic; then $A=C(X)\rtimes_\af \Z$ has tracial rank zero
      if the rotation numbers of $(X, \af)$ contain an irrational value (\cite{LP2}).
      Consequently every irrational rotation algebra has tracial rank zero. Thus, together with Theorem \ref{Ltr=0},
      we recover the result of Evans-Elliott (\cite{EE}) that every irrational rotation \CA\, is an inductive limit of circle algebras
      which has real rank zero.
     Theorem \ref{TLP2} brought a large part of crossed product \CA s into the classifiable class of simple \CA s.

    Later, with
    Theorem \ref{Ltr=0},
    we move   beyond the class  of  simple \CA s of real rank zero.

      \begin{thm}[Lin \cite{LinTAI}]\label{Ltr=1}
      Let  $A$ and $B$ be unital simple \CA s of tracial rank at most one in the UCT class.
      Then $A\cong B$ if and only if
      \beq\nonumber
       {\rm Ell}(A)\cong {\rm Ell}(B).
       \eneq
      \end{thm}
      It was proved in \cite{G5} that a unital simple AH-algebra with  very slow dimension growth
      has tracial rank at most one.
  The classification of C*-algebras which can be tracially approximated by splitting interval algebras was also obtained in \cite{Niu-thesis}, \cite{Niu-TAS-I}, \cite{Niu-TAS-II}, and \cite{Niu-TAS-III}.

      Applying a remarkable work of Winter's, his deformation result  (\cite{Wloc}), as pointed out, by Winter,
      led to the classification of ${\cal Z}$-stable \CA s  $A$ with the property
      that $A\otimes U$ has finite tracial rank for an arbitrary infinite dimensional
      UHF-algebra $U$. 
      
       Let us say, from now on,  that 
      $A$ has rational tracial rank $k,$ if $A\otimes F$ has tracial rank $k$
      for some infinite dimensional UHF-algebra $F.$

\begin{thm}[Lin-Niu \cite{LN-1}]\label{TLN-1}
Let $A$ and $B$ be unital separable simple ${\cal Z}$-stable \CA s in the UCT class.
Suppose that $A$ and $B$ has rational tracial rank zero.
Then
$A\cong B$ if and only if
\beq
{\rm Ell}(A)\cong {\rm Ell}(B).
\eneq
\end{thm}

One of the important features of the theorem above is that
the class of $\mathcal Z$-stable unital separable simple \CA  s with rational tracial rank zero
includes the Jiang-Su algebra ${\cal Z}$ (i.e.~ $TR(\mathcal Z\otimes F) = 0$ if $F$ is UHF), and it is the only algebra in the class with that invariant (the same as $\C$).
%

Toms and Winter proved the following generalization of Theorem \ref{TLP2} above.

\begin{thm}[Toms-Winter, \cite{TW1}]\label{TW1}
Let $X$ be an infinite, compact, finite dimensional metrizable space,   and
      $\af: X\to X$ be a minimal homeomorphism and
     $A=C(X)\rtimes_\af\Z.$
     Then $A$ has rational tracial rank zero if and only
     if its projections separate traces.
\end{thm}

One may note that Theorem \cite{TW1}  above brought these crossed products into
the classifiable \CA s, in particular, all crossed products \CA s of the form
$C(X)\rtimes_\af\Z$ for which $\af$ is minumal and uniquely ergodic.

With further development of the so-called Basic Homotopy Lemma
(\cite{Linhomtopytr=1})
and asymptotic unitary equivalence of \hm s (see \cite{Lnamj}), we reached the following
classification result.

\begin{thm}[Lin \cite{Lininv}]\label{Lininv}
Let $A$ and $B$ be unital separable amenable simple ${\cal Z}$-stable
\CA s with rational tracial rank at most one in the UCT class.
Then $A\cong B$ if and only if
\beq
{\rm Ell}(A)\cong {\rm Ell}(B).
\eneq
\end{thm}

We also have the following range theorem.

\begin{thm}[Lin-Niu \cite{LN-2}]\label{L-N-2}
Let $(G_0, (G_0)_+, u)$ be a countable ordered weakly unperforated and rationally Riesz group  (see  Proposition 5.8  and Definition 5.4 of \cite{LN-2}), let $G_1$ be a countable abelian group, let $\Delta$ be a metrizable Choquet simplex and let
$\rho: G_0\to \Aff(\Delta)$ be an order preserving \hm\, inducing a surjective affine continuous map
$\rho^\sharp: \Delta \to S(G_0, u)$ such that $\rho^\sharp(\partial_e(\Delta))=\partial_e(S(G, u)).$
Then there exists one (and exactly one, up to isomorphism) unital simple amenable $\mathcal Z$-stable C*-algebra $A$ with rational tracial rank at most one such that $$\mathrm{Ell}(A) = ((G_0, (G_0)_+, u), G_1, \Delta, \rho).$$
\end{thm}

This range theorem confirms that the class of rational tracial rank at most one
contains vast members of ${\cal Z}$-stable \CA s.

Unfortunately, this range theorem does not cover all possible Elliott invariants for
unital separable simple ${\cal Z}$-stable \CA s.  There are examples of unital separable simple
${\cal Z}$-stable \CA s $A$
for which  $K_0(A)$ does not have the rational Riesz property  and the pairing $\rho^\sharp$ does not map
extreme points to extreme points.


Let us now mention another non-commuative dimension for \CA s.


\begin{df}
Let $A$ and $B$ be C*-algebras. Recall (\cite{KW}) that a completely positive map $\phi: A\to B$ is said to have order zero if
$$ab=0\  \Rightarrow\  \phi(a)\phi(b)=0,\quad a, b\in A.$$
\end{df}

\begin{df}[\cite{KW}, \cite{WZ-ndim}]\label{DefDr}
A C*-algebra $A$ has nuclear dimension at most $n$, if there exists a net $(F_\lambda , \psi_\lambda , \phi_\lambda)$, $\lambda\in\Lambda$, such that the $F_\lambda$ are finite dimensional C*-algebras, and such that $\psi_\lambda: A\to F_\lambda$ and $\phi_\lambda: F_\lambda \to A$ are completely positive maps satisfying \\
$~\ \ $(1) $\phi_\lambda\circ \psi_\lambda \to \id_A$ pointwise (in norm),\\
$~\ \ $(2) $\|\psi_\lambda\|\leq 1$, and\\
$~\ \ $(3) for each $\lambda$, there is a decomposition
$F_\lambda=F_\lambda^{(0)}\oplus\cdots\oplus F_\lambda^{(n)}$
such that each restriction
$\phi_\lambda|_{F_\lambda^{(j)}}$ is a contractive order zero map.

Moreover, if the the map $\phi_\lambda$ can be chosen to be contractive itself, then $A$ is said to have  decomposition rank at most $n$.

\end{df}

\begin{thm}[Winter \cite{Winter-Z-stable-02}]\label{Wn}
Let A be a separable, simple, nonelementary, unital \CA\,
with finite nuclear dimension. Then $A$ is ${\cal Z}$-stable.
\end{thm}

We will present a classification theorem for
separable amenable simple ${\cal Z}$-stable \CA s.

      \section{The end of the beginning, unital case}

Recall that  the Jiang-Su algebra ${\cal Z}$ (\cite{JS} and see Theorem \ref{JS} above) constructed  in 1998 during the development of the Elliott program is an
infinite dimensional unital amenable  \CA\, in the UCT class with the feature
that its Elliott invariant
 is exactly the same as that of the complex field.
 As  stated earlier
 (see \cite{GJS} and Theorem \ref{GJS} above), for any separable simple \CA\,
$A,$ $A\otimes {\cal Z}$
 and $A$ have the same tracial structure and
$K_i(A)=K_i(A\otimes {\cal Z})$ ($i=0,1$). Moreover, $A\otimes {\cal Z}$ and $A$ are in the same  one of the three cases
mentioned at the end of Section 1.
 Since $A\otimes {\cal Z}$ is ${\cal Z}$-stable, i.e.,
$A\otimes {\cal Z}\cong (A\otimes {\cal Z})\otimes {\cal Z},$
naturally one studies  simple ${\cal Z}$-stable \CA s (see also \cite{Tex} for  further non-${\cal Z}$-stable examples).  Fortunately,
a separable nuclear  simple \CA\, is ${\cal Z}$-stable if and only if it has finite nuclear dimension (see
\cite{Winter-Z-stable-02}, \cite{MS},
\cite{MS2},
 \cite{SWW}, \cite{CETWW}, \cite{T-0-Z},  and \cite{CE}).

\begin{df}[See \cite{ET-PL} and  \cite{point-line}]\label{DfC1}
{\rm
Let $F_1$ and $F_2$ be two finite dimensional \CA s.
Suppose that there are two unital \hm s
$\phi_0, \phi_1: F_1\to F_2.$
Consider  the mapping torus \index{$A(F_1, F_2,\phi_0, \phi_1)$}  $M_{\phi_0, \phi_1}:$
$$
A=A(F_1, F_2,\phi_0, \phi_1)
=\{(f,g)\in  C([0,1], F_2)\oplus F_1: f(0)=\phi_0(g)\andeqn f(1)=\phi_1(g)\}.
$$
These \CA s {{were}} introduced into the Elliott program by Elliott and Thomsen (\cite{ET-PL}), and in \cite{point-line}, Elliott used this class of \CA s {{and some other building blocks with 2-dimensional spectra}} to realize any weakly unperforated {{simple}} ordered group {{with order unit}} as the $K_0$-group of a simple ASH \CA. Denote by ${\cal C}$ the class of all unital \CA s of the form $A=A(F_1, F_2, \phi_0, \phi_1)$ (which includes  all finite dimensional \CA s).
These \CA s will be called Elliott-Thomsen building blocks as well as
1-dimensional noncommutative CW complexes.

Let $\lambda: A\to C([0,1], F_2)$ be defined by $\lambda((f,a))=f.$
There is a canonical map $\pi_e: A \to F_1$ defined by $\pi_e(f,g)=g$ for all
pair $(f, g)\in A.$  It is a surjective map.\index{$\pi_e$}

If $A\in {\cal C}$, then $A$ is the pull-back  {{corresponding to the diagram}}
\begin{equation}\label{pull-back}
\xymatrix{
A \ar@{-->}[rr]^{\lambda} \ar@{-->}[d]^-{\pi_e}  && C([0,1], F_2) \ar[d]^-{(\pi_0, \pi_1)} \\
F_1 \ar[rr]^-{(\phi_0, \phi_1)} & & F_2 \oplus F_2\,.
}
\end{equation}

}
\end{df}

\begin{df}\label{Dgtrf}
Let $A$ be a  unital separable simple \CA. We say $A$ has generalized tracial rank at most one,
if, for any $\ep>0,$ any finite subset ${\cal F}\subset A$ and any $a\in A_+\setminus \{0\},$
there is a \SCA\, $B\in {\cal C}$ such that $p=1_B$ and

(i) $\|px-xp\|<\ep\rforal x\in {\cal F},$

(ii) $pxp\in_\ep B\rforal x\in {\cal F}$, and

(iii) $p\lesssim a.$

If $A$ has generalized tracial rank at most one, we may write
$gTR(A)\le 1.$  If there exists an infinite dimensional UHF-algebra
$F$ such that $gTR(A\otimes F)\le 1,$ we say $A$ has rationally generalized tracial rank at most one.

\end{df}

The introduction of the class of unital separable amenable simple \CA s of rational
generalized tracial rank at most one
is important.  From the following range theorem, we conclude that
this class of simple \CA s exhausts all possible Elliott invariants for unital separable simple
${\cal Z}$-stable \CA s.

\begin{thm}[Gong-Lin-Niu Theorem 13.50 of \cite{GLN1}]\label{range 0.35}
For any simple weakly unperforated Elliott invariant $\big((G,G_+,u),~
K, \DT, r \big)$, where $(G, G_+, u)$ is a weakly unperforated countable simple ordered
group $G$ with positive cone $G_+$ and with a non-zero element $u,$ $K$ is any
countable abelian group, $\DT$ any compact metrizable Choquet simplex, and
$r: G\to \Aff(\DT)$ any order preserving \hm\, (with the strict order on $\Aff(\DT)$),
 there is a
 unital simple ${\cal Z}$-stable \CA\, $A$  with rational tracial rank at most one
which is also an
$ASH$-algebra  such that
$$\big((K_0(A),K_0(A)_+, [1_A]),~ K_1(A), T(A), \rho_A \big)\cong
\big((G,G_+,u),~ K, \DT, r \big).$$
\end{thm}

\begin{rem}
In the case of finite separable simple unital ${\cal Z}$-stable \CA s,
$K_0(A)_+$ is determined by $[1_A]$ and the pairing $\rho_A: K_0(A)\to \Aff(T(A).$
Therefore the invariant could be simplified as
\beq
(K_0(A), [1_A], K_1(A), T(A), \rho_{A}).
\eneq
However, we keep the positive cone $K_0(A)_+$ following the early tradition.
\end{rem}

Let us present the following isomorphism theorem:

\begin{thm}[Gong-Lin-Niu \cite{GLN2}]\label{TGLN}
Let $A$ and $B$ be unital separable amenable simple ${\cal Z}$-stable \CA s with rational generalized tracial rank at most one
in the UCT class.  Then
$A\cong B$ if and only if
\beq
{\rm Ell}(A)\cong {\rm Ell}(B).
\eneq
\end{thm}

The question is then whether  every finite unital separable amenable simple ${\cal Z}$-stable \CA\,
has rational generalized tracial rank at most one.

%

\begin{thm}
[Castillejos, Evington, Tikuisis, White, and Winter \cite{CETWW}]\label{CETWW}
Let A be an infinite dimensional, separable, simple, unital,
nuclear \CA. Then the following statements are equivalent:

(i) $A$ has finite nuclear dimension;

(ii) $A$ is ${\cal Z}$-stable, i.e.,  $A\cong  A\otimes {\cal Z}.$
\end{thm}


Let us state the following theorem.

\begin{thm}[Elliott-Gong-Lin-Niu  Theorem 4.9 of \cite{EGLN-3}]\label{T2}
Let $A$ be a unital  simple separable amenable (non-zero) C*-algebra  which satisfies the UCT.
Then the following  properties
are equivalent:

{\rm (1)} ${\mathrm{gTR}}(A\otimes Q)\le 1;$

{\rm (2)} $A\otimes Q$ has finite nuclear dimension and ${\rm T}(A\otimes Q)={\rm T}_{{\rm qd}}(A\otimes Q)\not=\O;$

{\rm (3)} the decomposition rank of
$A\otimes Q$
is finite.

\noindent (Here $Q$ is the universal UHF algebra and ${\rm T}_{{\rm qd}}(A\otimes Q)$ is the set of quasidiagonal traces.)
\end{thm}

It follows from a result of
Tikuisis, White, and  Winter
\cite{TWW} that
every tracial state on a unital simple separable amenable
\CA\, which satisfies the UCT is quasidiagonal.   Since $A\otimes Q$ is always
${\cal Z}$-stable, by Theorem \ref{CETWW}, $A\otimes Q$  has finite nuclear dimension.
Therefore (2) holds for any finite unital separable amenable simple \CA\, $A$ in the UCT class, whence 
(1) in Theorem \ref{T2} always holds.
From this 
 we obtain  the following theorem.

\begin{thm}[Elliott-Gong-Lin-Niu \cite{EGLN-3}]\label{TEGLN-1}
Let $A$ be a unital simple separable amenable
\CA\, 
which satisfies the UCT.
Then  $A$ has rational  generalized tracial rank at most one.
%
\end{thm}


 Together with the isomorphism theorem
\ref{TGLN}, as well as the classification of purely infinite simple \CA s  of Kirchberg and
Phillips \cite{KP} and \cite{Pclass}),  we further obtain the following:

\begin{cor}[Elliott-Gong-Lin-Niu \cite{EGLN-3}]\label{EGLN3}
Let $A$ and $B$ be unital separable simple amenable ${\cal Z}$-stable \CA s in the UCT class.
Then
$A\cong B$ if and only if
\beq
{\rm Ell}(A)\cong {\rm Ell}(B).
\eneq
\end{cor}





\section{The end of the beginning, non-unital case}

Next we will discuss the non-unital case.
We will mainly discuss the stably projectionless case.
In fact, the classification  of non-unital separable amenable simple ${\cal Z}$-stable
\CA s in the UCT class which are not stably projectionless can easily  be reduced to the unital
case (see  the proof of Theorem \ref{Tnonunitalproj} below).

However, there are separable simple \CA s  $A$ which are
inductive limits of sub-homogeneous \CA s such that $A\otimes {\cal K}$ admits no non-zero
projections.  One of them, $W,$ a Razak algebra, has a unique tracial state and
has zero $K_0$ and zero $K_1$ (\cite{RzW}).  In other words, Case 2 mentioned in Section 1 does
occur.

{{The study of stably projectionless simple \CA s actually  has a long history (see  {{Kishimoto's 1980 paper \cite{K80},
and Blackadar and Cuntz's 1982 paper \cite{BC},}} for example).}}
Stably projectionless simple \CA s can naturally {{occur, for example,}} in the study of
flow actions (see \cite{K80}, \cite{KK}, 
\cite{Rl}, and \cite{T-0-Z}).
Tsang (\cite{Tsang}) showed that any metrizable Choquet simplex
can be the tracial state space for  some  stably projectionless simple  amenable  \CA.
In fact, in \cite{GLII}, we show that there is a unique
separable stably projectionless simple \CA\, $\zo$ with finite nuclear dimension in the UCT class
and with $K_0(\zo)=\Z,$ $K_1(\zo)=\{0\}$, and with a unique tracial state (such \CA s {{were}} known to exist).
It turns out that, for any separable simple  nuclear \CA\, $A,$ $A\otimes \zo$ is stably projectionless.}}
Furthermore,
for any abelian group $G_0,$ any compact metrizable Choquet simplex $\Delta$, and
any \hm\, $\rho: G_0\to \Aff(\Delta),$ the space of continuous real affine functions on $\Delta$, such that $\rho(G_0)\cap \Aff_+(\Delta)=\{0\},$  there exists a  separable  stably projectionless simple nuclear  ${\cal Z}$-stable \CA\, $A$
such that $(K_0(A), T(A), \rho_A)=(G_0, \Delta, \rho)$ (and with arbitrarily given $K_1(A)$)  {{(see \cite{point-line} and \cite{GLrange}).}}

Since the stably projectionless  case is not included in \cite{GLN1} and \cite{GLN2}, we will provide a few more details.


\subsection{The Elliott invariant and its range}

We start with the description of the invariant.

\begin{df}\label{DTtilde}
Let $A$ be a \CA.  Denote by $T(A)$ the tracial state {{space}} of $A$ {{(which could be an empty set).}}
{{Denote by $T_f(A)$ the set of faithful traces of $A$
(if $a\in A_+\setminus \{0\}$ and
$\tau\in T_f(A),$ then $\tau({{a}})>0$).}}
Let ${\tilde{T}}(A)$ denote the cone of densely defined,
positive lower semi-continuous traces on $A$ equipped with the topology
of point-wise convergence on elements of the Pedersen ideal  ${\rm Ped}(A)$---the minimal dense hereditary ideal of $A.$
Let $B$ be another \CA\, with $\td T(B)\not=\{0\}$
and let $\phi: A\to B$ be a \hm.

{\it{In what follows we will also write $\phi$ for $\phi\otimes \id_{M_k}: M_k(A)\to M_k(B)$
whenever it is convenient.}}
We will write   $\phi_T: \td T(B)\to \td T(A)$ for the induced continuous affine map.
Denote by $\td T^{b}(A)$ the subset of $\td T(A)$ of traces which are bounded on $A.$
Of course $T(A)\subset \td T^b(A).$ Set $T_0(A):=\{\tau\in \td T(A): \|\tau\|\le 1\}.$ It is a compact convex subset of $\td T(A).$

Let $r\ge 1$ be an integer and $\tau\in {\tilde T}(A).$
We will continue to {{write}} $\tau$  on $A\otimes M_r$ for $\tau\otimes {\rm Tr},$ where ${\rm Tr}$ is the standard
trace on $M_r.$
Let  $S$ be a convex subset (of a convex topological cone $C$ with a Choquet simplex as a base).
{{We assume  that the convex cone contains 0. Denote by
$\Aff(C)$ the set of affine real continuous functions on $C$ which vanish at $0.$}}
Define (see \cite{Rl})
\beq
\Aff(S)&=&\{f|_S: f\in \Aff(C)\}\andeqn \Aff(S)_+=\{f\in  \Aff(S): f\ge 0\},\\
\Aff_+(S)&=&\{f\in \Aff(S):
 f(\tau)>0\,\,{\rm for}\,\,\tau\in S\setminus \{0\}\}\cup \{0\},\\
{\rm LAff}(S)_+&=&\{f:S\to [0,\infty]: \exists \{f_n\}, f_n\nearrow f,\,\,
 f_n\in \Aff(S)_+\},\\
{\rm LAff}_+(S)&=&\{f:S\to [0,\infty]: \exists \{f_n\}, f_n\nearrow f,\,\,
 f_n\in \Aff_+(S)\}.
 \eneq
 To include the case that $S=\{0\},$ we view that it is the only point in the zero cone.
 In that case, we define  $\Aff(S)=\Aff_+(S)=\{0\}={\rm LAff}_+(S).$
 For the greater part of this paper, $S={\widetilde T}(A)$ {{or}} $S=T(A),$ or $S=T_0(A)$
in the above definition will be used.
Recall {{that}} $0\in {\tilde T}(A)$ and if ${{g}}\in \LAff(\td T(A)),$ then ${{g}}(0)=0.$


\end{df}

%
%
\begin{df}\label{Dfep}
For any $\ep>0,$ define $f_\ep\in {{C([0,\infty))_+}}$ by
$f_\ep(t)=0$ if $t\in [0, \ep/2],$ $f_\ep(t)=1$ if $t\in [\ep, \infty)$, and
$f_\ep(t)$ is linear in $(\ep/2, \ep).$

Let $A$ be a \CA\, and $\tau$ a trace
on $A$.
For each $a\in A_+,$ define $d_\tau(a)=\lim_{\ep\to 0} {{\tau(f_\ep(a)).}}$
Note that $f_\ep(a)\in {\rm Ped}(A)$ for all $a\in A_+$ {{(see Theorem 5.6.1 of  \cite{Pbook} and its proof)}}.

{{Let $S$ be a  convex subset of $\widetilde{T}(A)$ and $a\in M_n(A)_+.$
The function $\hat{a}(s)=s(a)$ (for $s\in S$) is an  affine function.
Define $\widehat{\la a\ra }(s)=d_s(a)=\lim_{\ep\to 0}s(f_\ep(a))$ (for $s\in S$),
which is a lower semicontinuous function.
If $a\in {\rm Ped}(A)_+,$ then $\hat{a}$ is in $\Aff(S)_+$ and
$\widehat{\la a\ra}\in \LAff(S)_+$  {{in general}} (see \ref{DTtilde}). 
If $A$ is simple, then $\widehat{\la a\ra}\in \LAff_+(S).$ {{Note
that $\hat{a}$ is different from $\widehat{\la a\ra}.$}}
 In most cases, $S$ is ${\tilde T}(A),$ $T_0(A),$ or $T(A).$}}
Note {{also}} that, when $A$ is simple,
 there is a natural map from ${\rm Cu}(A)$ to {{${\rm LAff}_+(\tilde{T}(A))$}}, sending $\la a \ra$ to $\widehat{\la a \ra}$.

\end{df}

\begin{NN}\label{range4.1}
If $A$ is a unital \CA\, and $T(A)\not=\emptyset,$ then
there is a canonical  
\hm\, $\rho_A: K_0(A)\to \Aff(T(A)).$
%
{{Now consider the case that $A$ is not unital.
Let $\pi_\C^A: \td A\to \C$ denote the quotient map.
Suppose that $T(A)\not=\emptyset.$
Let $\tau_\C:=\tau_\C^A: \td A\to \C$ denote
the tracial state which factors through $\pi_\C^A.$
Then
\beq
T(\td A) = \big\{t{{\tau_{_\C}^A}}+(1-t)\tau:~ t\in [0,1],~ \tau\in {{T(A)}}\big\}.
\eneq
The map ${{T(A)}}\hookrightarrow {{T(\td A)}}$ induces a map {{${{\Aff(T(\tilde{A}))}}\to \Aff(T(A))$}}. Then
 the  map {{$\rho_{\td A}: K_0({{\td A}})\to \Aff(T(\td A))$}} induces a \hm\,  $\rho':~ K_0(A)\to \Aff (T(A))$ by
\beq
{{\rho':~K_0(A)\to  K_0({{\td A}})\stackrel{\rho_{\td A}}{\longrightarrow}
\Aff(T(\td A))   \to \Aff(T(A)).}}
\eneq
However, in the case that $A\not={\rm Ped}(A),$ we will not use $\rho'$  in general, as it is possible
that $T(A)=\emptyset$ but ${\td T}(A)$ is rich (consider the case $A\cong A\otimes {\cal K}$).}}
\end{NN}

\begin{df}{{[Definition 2.15 of \cite{GLrange}]}}\label{Dparing}
{{Let $A$ be a \CA\, with $\td T(A)\not=\{0\}.$
If $\tau\in \td T(A)$ is bounded on $A,$ then $\tau$ can be extended naturally
to {{a trace}} on $\td A.$
Recall that $\td T^b(A)$ is the set of bounded traces on $A.$
Denote by $\rho_A^b: K_0(A)\to \Aff(\td T^b(A))$ the \hm\, defined by
$\rho_A^b([p]-[q])=\tau(p)-\tau(q)$ for all $\tau\in \td T^b(A)$
and for projections $p, q\in M_n(\td A)$ (for some integer $n\ge 1$)
and $\pi_\C^A(p)=\pi_\C^A(q).$ {{Note $p-q\in  M_n(A).$
Therefore $\rho_A^b([p]-[q])$ is continuous on $\td T^b(A).$}}
In the case that $\td T^b(A)=\td T(A),$ for example, $A={\rm Ped}(A),$
we write $\rho_A:=\rho_A^b.$

Let $A$ be a $\sigma$-unital \CA\, with a strictly positive element $0\le e\le 1.$
Put  $e_n:=f_{1/2^n}(e).$ Then $\{e_n\}$ forms an approximate identity
for $A.$ Note $e_n\in {\rm Ped}(A)$ for all $n.$
Set $A_n={\rm Her}(e_n):=\overline{e_nAe_n}.$
Denote by $\iota_n:A_n\to A_{n+1}$   and $j_n: A_n\to A$ the embeddings. These extend to
$\iota_n^\sim: \td A_n\to \td A_{n+1}$ and $j_n^\sim :\td A_n\to \td A$ {{unitally.}}
Note that $e_n\in {\rm Ped}(A_{n+1}).$ Thus  $\iota_n$  and $j_n$ induce
continuous cone maps ${\iota_n}_T^b: \td T^b(A_{n+1})\to
\td T^b(A_n)$ and ${j_n}_T: \td T(A)\to \td T^b(A_n)$
(defined by ${\iota_n}_T^b(\tau)(a)=\tau(\iota_n(a))$
for $\tau\in {{\td T^b}}(A_{n+1}),$ and ${j_n}_T(\tau)(a)=\tau(j_n(a))$
for all $\tau\in \td T(A)$ and all $a\in A_n$), respectively.
Denote by $\iota_n^\sharp: \Aff(\td T^b(A_n))
\to {{\Aff(\td T^b(A_{n+1}))}}$  and $j_n^\sharp: \Aff(\td T^b(A_n))\to \Aff(\td T(A))$ the induced continuous linear maps.
Recall that $\bigcup_{n=1}^\infty A_n$ is dense in ${\rm Ped}(A).$
A direct computation shows that one may obtain  the following inverse direct limit of {{convex  topological cones}} (with continuous cone maps):
\beq\label{Drho-trace}
\td T^b(A_1)\stackrel{{\iota_1}_T^b}
{\longleftarrow} \td T^b(A_2)\stackrel{{\iota_2}_{T}^b}{\longleftarrow} \td T^b(A_3)
\cd\longleftarrow\cd\longleftarrow \td T(A),
\eneq
}}
which induces the following commutative {{diagram:}}
\beq
\Aff(\td T^b(A_1))\stackrel{\iota_1^\sharp}
{\longrightarrow} \Aff (\td T^b(A_2))\stackrel{{\iota_2}^\sharp}{\longrightarrow} \Aff (\td T^b(A_3))
\cd\longrightarrow\cd\longrightarrow \Aff (\td T(A)).
\eneq
Hence one also has the following commutative diagram:
\begin{displaymath}
    \xymatrix{
        K_0(A_1) \ar
        [d]_{\rho_{A_1}}\ar[r]^{\iota_{1*o}} & K_0(A_2) \ar[r]^{\iota_{2*o}} \ar
        [d]_{\rho_{A_2}}& K_0(A_3) \ar[r]
        \ar
        [d]_{\rho_{A_3}}& \cd K_0({{A}} )\\
        \Aff(\td T^b({{A}}_1))
        \ar[r]^{\iota^{\sharp}_{1,2}}
        &
         \Aff(\td T^b({{A}}_2)) \ar[r]^{\iota_2^{\sharp}}
         &
         \Aff(\td T^b({{A}}_3))
         \ar[r]
         & \cd
         \Aff(\td T({{A}})).}
\end{displaymath}
Thus one obtains a \hm\, $\rho: K_0(A)\to \Aff(\td T(A)).$
It should be noted that, when $A$ is simple, $\td T^b(A_n)=\td T(A_n)$ for all $n$
{{(see Definition {{2.15}} of \cite{GLrange} for more details). Moreover, the map $\rho$ does not depend
on the choice of $\{e_n\}.$}} {{We will write $\rho_A:=\rho.$ In  the case that $T(A)$ generates ${\tilde T}(A)$
such as the case that {{$A$ has continuous scale,}}
we may also write $\rho_A: K_0(A)\to \Aff(T(A))$  by   restricting
$\rho_A(x)$  to $T(A)\subset {\tilde T}(A)$ and for all $x\in K_0(A).$}}
\end{df}


Now let us present an isomorphism theorem
for non-unital separable amenable simple ${\cal  Z}$-stable \CA s in the UCT class
which are not stably projectionless.
Suppose that $A$ has a strictly positive element $e_A.$
Denote by $\widetilde{\la e_A\ra}$ the dimension function
$\widetilde{\la e_A\ra}(\tau)=d_\tau(\la e_A\ra)$ ($\tau\in \widetilde{T}(A)$).
Note that the dimension function $\widetilde{\la e_A\ra}$ is independent of the choice
of the strictly positive element $e_A.$

\begin{thm}\label{Tnonunitalproj}
Let $A$ and $B$ be non-unital separable simple ${\cal Z}$-stable \CA s in the UCT class with  strictly positive elements $e_A$ and $e_B,$ respectively. Suppose that $A\otimes {\cal K}$ contains a non-zero projection.
Then $A\cong B$ if and only if
\beq\label{elliso-1}
(K_0(A), \widetilde{\la e_A\ra}, \widetilde{T}(A), \rho_A, K_1(A))\cong
(K_0(B), \widetilde{\la e_B\ra}, \widetilde{T}(B), \rho_B, K_1(B)).
\eneq
\end{thm}

\begin{df}

Note that \eqref{elliso-1} means that there is an
a group isomorphism  $\kappa_0: K_0(A)\to K_0(B),$
an affine homeomorphism $\kappa_T: \widetilde{T}(A)\to \widetilde{T}(B)$
such that
\beq
\rho_B(\kappa_0(x))(\tau)=\rho_A(x)(\kappa_T^{-1}(\tau))\rforal x\in K_0(A)\andeqn \tau\in \widetilde{T}(B),
\eneq
and $\widetilde{\la e_A\ra}(\kappa_T^{-1}(\tau))=\widetilde{\la e_B\ra}(\tau)$ for all $\tau\in \widetilde{T}(B),$ and a group isomorphism
$\kappa_1: K_1(A)\to K_1(B).$
\end{df}

{\bf Proof of Theorem \ref{Tnonunitalproj}: }
Let
\beq\nonumber
\Gamma: (K_0(A), \widetilde{\la e_A\ra}, \widetilde{T}(A), \rho_A, K_1(A))\cong
(K_0(B), \widetilde{\la e_B\ra}, \widetilde{T}(B), \rho_B, K_1(B))
\eneq
be an isomorphism.

Let $p\in A\otimes {\cal K}$ be a non-zero projection and set
$A_1=p(A\otimes {\cal K})p.$ Then $A_1$ is a unital simple separable amenable 
${\cal Z}$-stable \CA\, in the UCT class.
Let $q\in B\otimes {\cal K}$ be a projection
such that $[q]=\kappa_0([p]).$  Put $B_1=q(B\otimes {\cal K})q.$
By applying Theorem \ref{EGLN3},   there  is an isomorphism
$\phi: A_1\to B_1$ such that
$\phi$ induces  $\Gamma.$
Note that $A$ is non-unital, and since
${\rm Ell}(B)\cong {\rm Ell}(A),$ $\Sigma(K_0(B))$ has either no maximum element
or the maximum element $u\in \Sigma(K_0(B))$ does not have the property
$\rho_B(u)=\widetilde{\la e_B\ra}.$ Hence $B$ is not unital.  Consequently
neither $\la e_A\ra$ nor $\la e_B\ra$ can be represented by projections.
We may view $e_A\in A_1\otimes {\cal K}$ and $e_B\in B_1\otimes {\cal K}$
with 
$\overline{e_A(A_1\otimes {\cal K})e_A}=A$ and
$\overline{e_B(B_1\otimes {\cal K})e_B}=B.$
On the other hand, we must have  $A\cong \overline{\psi(e_A)(B_1\otimes {\cal K})\psi(e_A)}$ and
$\widetilde{\la \psi(e_A)\ra}=\widetilde{\la e_B\ra},$ and both $\psi(e_A)$ and $\la e_B\ra$
are not represented by projections. Then, by  Theorem 6.6 of \cite{ESR-Cuntz},
$\la \psi(e_A)\ra=\la e_B\ra.$   It follows from Proposition 3.3 of \cite{Rlz} (see also Theorem 3 of \cite{CES}
as well as  Proposition 4.3 of \cite{ORT})
that there exists a partial isometry $v\in (B_1\otimes {\cal K})^{**}$ such that $v
\overline{\psi(e_A)(B_1\otimes {\cal K})}\in B_1\otimes {\cal K},$ $v^*\overline{e_B(B_1\otimes {\cal K})e_B}
\in B_1\otimes {\cal K},$ and
$v^*v=p,$ the open projection corresponding to $\psi(e_A)$, and
$vv^*=q,$ the open projection corresponding to $e_B.$
Then $\phi_1: A\to B$ defined by
$\phi_1(a)=v\psi(a) v^*$ for all $a\in A$ is an isomorphism from
$A$ to $B$ which carries $\Gamma.$
%


\vspace{0.2in}

With Corollary \ref{EGLN3} and Theorem \ref{Tnonunitalproj},  we complete the classification
of Case (1).
The world is not perfect. As mentioned earlier, Case 2 does occur.
Let us now discuss the Case 2, the stably projectionless case.

\begin{NN}\label{DElliott+}

{{Recall from Definition {{2.7}} of \cite{GLrange} {{(see (2.8) of \cite{GLrange} also)}}, {{that}} a (simple)}} ordered group pairing is a triple $(G, T,  \rho),$ where
$G$ is a countable abelian group, $T$ is  {{a convex  topological cone}} with a Choquet simplex as
base,
and  $\rho: G\to \Aff(T)$  is a \hm.
Let us also add the case that $T=\{0\}$ with $\rho=0.$
{{Define
\beq\label{++}
\Aff_+(T)^+=\{f\in \Aff_+(T): f(t)>0, \,\,\, {\rm whenever}\,\,\, t\not=0\}.
\eneq
Note that, in the case that $T=\{0\},$ $\Aff_+(T)^+=\Aff_+(T)=\{0\}.$}} 

Define $G^+=\{g\in G: \rho(g)\in \Aff_+(T)^+\}\cup \{0\}.$ When $T\not=\{0\},$ if $G^+\not=\{0\},$
then $(G, G^+)$ is an ordered group {{in the sense that $G^+\cap (-G^+)=\{0\}$ and $G^+-G^+=G$}}. It has the property that if $ng>0$ for some integer $n>0,$
then $g>0.$  In other words, $(G, G^+)$ is weakly unperforated.
When $T=\{0\},$ recall that $\Aff_+(T)^+=\{0\},$  and hence  $G^+=G.$

The ordered group pairing {{mentioned}} above satisfies the following condition:  when $G\not=G^+,$
either
 $G^+=\{0\},$ or
$(G, G^+)$ is a simple ordered group, i.e.,  every element $g\in G^+\setminus \{0\}$
is  an order unit.
%

{\em A simple ordered group pairing with a distinguished set,
or a scaled simple  ordered group pairing } is a quintuple $(G, \Sigma_\rho(G), T, s, \rho)$ such that
$(G, T, \rho)$ is a simple ordered group pairing,
 $s\in \LAff_+(T)\setminus \{0\}$, or, $s=0,$ when $T=\{0\},$  and 
\beq
\Sigma_\rho(G):=\{g\in G^+: \rho(g)<s\}\,\,\,
{\rm or}\,\,\,
\Sigma_\rho(G):=\{g\in G^+: \rho(g)<s\}\cup\{u\},
\eneq
where $\rho(u)=s$ and  $\rho(g)<s$ means $s-\rho(g)\in \Aff_+(T)$ and $s-\rho(g)\not= 0.$
Here $\Sigma_\rho(G)$ is the distinguished set. It could contain only one element 
or even be an empty set.
We also allow $\Sigma_\rho(G)=\{0\}.$     
In the case that $T=\{0\},$ $s=0.$
No element $g\in G^+$ with $\rho(g)<s.$ Hence 
when $T=\{0\},$
if 
$\Sigma_\rho(G)=\{g\in G^+: \rho(g)<s\},$ then  $\Sigma_\rho(G)=\emptyset,$ and
if $\Sigma_\rho(G)=\{g\in G^+:\rho(g)<s\}\cup u$ 
with $\rho(u)=s,$ then $\Sigma_\rho(g)=\{u\}.$ Therefore, when $T=\{0\},$
there are two cases, either the distinguished set $\Sigma_\rho(G)=\emptyset,$ or $\Sigma_\rho(G)=\{u\},$ a single 
element (which could be zero).

Note {{that}} $s(\tau)$ could be infinite for some $\tau\in T.$
It is called a unital scaled simple ordered group pairing, if $\Sigma_\rho(G)=\{g\in G^+: \rho(g)<s\}\cup \{u\}$
with $\rho(u)=s,$ in which case, $u$ is called the unit for $G.$  Note {{also}} that, in this {{case,}}
$u$ is the maximum element of {{$\Sigma_\rho(G),$
one may}} write
$(G, u, T, \rho)$ for $(G, \Sigma_\rho(G), T, s,\rho).$  {{On the other hand,
$\Sigma_\rho(G)$ is determined by $s$ if it has no unit.}}
One may {{then}} write $(G, T, s, \rho)$ for $(G, \Sigma_\rho(G), T, s, \rho)$
(see Theorem {{5.2 of \cite{GLrange}.}}) 

Let $(G_i, \Sigma_\rho(G_i), T_i, s_i, \rho_i),$  $i=1,2,$ be {{scaled}} simple ordered group pairings.
A map
$$\Gamma_0: (G_1, \Sigma_\rho(G_1), T_1, s_1, \rho_1)\to (G_2, \Sigma_\rho(G_2), T_2, s_2, \rho_2)$$
is said to be an isomorphism, if it consists of
a group isomorphism
$\kappa_0: G_1\to G_2$ and an affine cone isomorphism $\kappa_T: T_1\to T_2$
such that
\beq\label{july12-2021}
{{\rho_2}}(\kappa_0(g))(t)={{\rho_1(g)}}(\kappa_T^{-1}(t))\rforal g\in G_1\andeqn t\in T_2,\andeqn\\ \label{july12-2021-1}
\kappa_0(\Sigma_\rho(G_1))=\Sigma_\rho(G_2)), \andeqn   s_1(\kappa_T^{-1}(t))=s_2(t) \rforal t\in T_2.
\eneq
\end{NN}

Recall that for
the classification of a class of \CA s, there must be {\em three} ingredients:

\noindent
(1) A definition of invariants for \CA s in the class; \\
 (2)  A (range) theorem  that  says exactly what values of the
invariants defined in (1)  occur for \CA s in the class;\\
 (3) the  isomorphism theorem, i.e., the statement that any two \CA s in the class with the same,
 i.e., isomorphic  invariants
are isomorphic.

  Let us now present a unified classification for finite separable amenable simple
  ${\cal Z}$-stable \CA s in the UCT class.

\begin{df}
From {{2.15 and 2.19}} of \cite{GLrange}, with the  pairing $\rho_A: K_0(A) \to {{\Aff (\td T(A))}}$ defined in \ref{Dparing} (see
2.15  of \cite{GLrange}),  the Elliott invariant  for separable simple \CA s {{(see {{\cite{Ellicm} and}} \cite{point-line}),}}
%
%
%
{{for}} the case  ${{{\tilde T}(A)}}\not=\{0\},$
 may be described by the following sextuple:
$$
{\rm{Ell}}(A):=((K_0(A), \Sigma_\rho(K_0(A)), {\tilde T}(A), \Sigma_A, \rho_A), K_1(A)),
$$
where $\Sigma_A$ is a function in $\LAff_+(\td T(A))$ defined by
\beq
\Sigma_A(\tau)=\sup\{\tau(a): a\in {\rm Ped}(A)_+,\,\,\|a\|\le 1\}.
\eneq
 Let $e_A\in A$ be a strictly positive element.
Then $\Sigma_A(\tau)=d_\tau(e_A)=\lim_{\ep\to 0} \tau(f_\ep(e_A))$ for all $\tau\in {{\td T(A),}}$ which
is independent of the choice of $e_A.$
In this case, if $\Sigma_\rho(K_0(A))$ has a maximum element $u$ and 
$\rho(u)=\Sigma_A,$ then,
\beq
\Sigma_\rho(K_0(A))=\{g\in K_0(A)^+: \rho(g)<\Sigma_A\}\cup\{u\}
\eneq
(recall that $K_0(A)^+=\{x\in K_0(A): \rho_A(x)\in \Aff_+(\widetilde{T}(A))^+\}\cup\{0\},$ see \eqref{++}).
Otherwise, 
\beq
\Sigma_\rho(K_0(A))=\{g\in K_0(A)^+: 
 \rho(g)<\Sigma_A\}.
\eneq
%
In the case that $\widetilde{T}(A)=\{0\}$ (hence $\rho_A=0$),  we also define 
\beq
{\rm Ell}(A)=\{K_0(A), \Sigma_\rho(K_0(A)), {\widetilde T}(A), \Sigma_A, \rho_A), K_1(A)),
\eneq
where $\Sigma_A=0.$  As in \ref{DElliott+}, since $\widetilde{T}(A)=\{0\},$  we have
$\{x\in K_0(A)^+: \rho(x)<0\}=\emptyset.$  Hence, in this case,
when $A$ has a unit,
\beq
\Sigma_\rho(K_0(A))=\{x\in K_0(A)^+: \rho(x)<0\}\cup\{[1_A]\}=\{[1_A]\},
\eneq
and, when $A$ is not unital,
\beq
\Sigma_\rho(K_0(A))=\{x\in K_0(A)^+: \rho(x)<0\}=\emptyset.
\eneq
Recall that the set $\Sigma_\rho(K_0(A))$ is the distinguished set.

If $A$ is a separable ${\cal Z}$-stable \CA\, with $\widetilde{T}(A)=\{0\},$ 
then by a result of M. R\o rdam  $A$ must be purely infinite (Corollary 5.4 of \cite{R4}). 
In this case,
\beq
&&{\rm Ell}(A)=(K_0(A), \{[1_A]\}, \{0\}, 0, 0, K_1(A)),\\
{\rm or}\,\,\, &&{\rm Ell}(A)=(K_0(A), \emptyset, \{0\}, 0, 0, K_1(A)).
\eneq
One may then simplify the invariant in these two cases by writing
\beq
{\rm Ell}(A)=(K_0(A), [1_A], K_1(A)),\,\,\,{\rm or}\,\,\, {\rm Ell}(A)=(K_0(A), \emptyset,  K_1(A)).
\eneq
In other words,  in the purely infinite simple case,
 $A$ has  no unit if and only if the distinguished set $\Sigma(K_0(A))=\emptyset.$
 The definition above gives the {\it first ingredient} of the classification.

Let $B$ be another separable  \CA.  A {{\hm\,}} $\Gamma: {\rm Ell}(A)\to {\rm Ell}(B)$ consists of a {{\hm\,}}   $\Gamma_0:(K_0(A), \Sigma_\rho(K_0(A)), {\widetilde T}(A), \Sigma_A, \rho_A)$ to $(K_0(B), \Sigma_\rho(K_0(B)), {\widetilde T}(B), \Sigma_B, \rho_B)$ (as
in \ref{DElliott+}, see (\ref{july12-2021}) and (\ref{july12-2021-1}) also) and a homomorphism  $\kappa_1:K_1(A)\to K_1(B)$. We say that $\Gamma$ is an isomorphism if  both $\Gamma_0$ and $\kappa_1$ are isomorphisms.

In the case that $\rho_A(K_0(A))\cap {\rm Aff}_+(\widetilde T(A))=\{0\},$ we often consider
the reduced case that $T(A)$ is compact which gives a base for $\widetilde T(A).$ In that  case, we
may write \\ ${\rm Ell}(A)=(K_0(A), T(A), \rho_A, K_1(A)).$ Note {{that}}, in the said situation,
$\Sigma_\rho(K_0(A))=\{0\}$ {{and}}
$\widetilde T(A)$ is determined by $T(A)$ and $\Sigma_A(\tau)=1$ for all $\tau\in T(A).$

\end{df}

Note, that in the special case  that $K_0(A)=K_1(A)=\{0\},$ the Elliott invariant 
may  be simplified to $(\widetilde{T}, \Sigma_A).$
Recall that a \CA\, $A$ is $KK$-contractible if it $KK$-equivalent to $0.$
In particular $K_i(A)=0$ ($i=0,1$).
Before we discuss the general case, let us first state the following theorem:

\begin{thm}[Elliott, Gong, Lin, and Niu, Theorem 7.5 of \cite{eglnkk0}]\label{TEGLNkk0}
The class of KK-contractible stably projectionless simple separable C*-algebras with finite nuclear
dimension is classified by the invariant  $(\widetilde{T}(A), \Sigma_A).$
Any \CA\, $A$ in this class is a simple inductive limit of Razak algebras.
\end{thm}
\begin{df}\label{dk0}
Let $A$ be a finite separable  is traditional to consider the following distinguished 
subset of $K_0(A),$ 
\beq
\Sigma(K_0(A))=\{x\in K_0(A): x=[p]\,\, {\rm for\,\, some\,\, projection},\,\, p\in A\}.
\eneq
Moreover $K_0(A)_+=\{x\in K_0(A): x=[p]\,\,{\rm for\,\, some\,\, projection}\}.$
\end{df}

Note that, for a finite exact separable simple ${\cal Z}$-stable \CA\, $A, $ by a result of 
 R\o rdam (Corollary. 5.4 of \cite{R4}), $\widetilde{T}(A)\not=\emptyset.$
For finite separable amenable simple \CA s, the next  theorem should be regarded as an integral part 
of the range theorem.

\begin{thm}[Gong-Lin, Theorem 5.2 of \cite{GLrange} and  Corollary A.7  of \cite{eglnp}]\label{GLrange1}
 Let $A$ be a finite separable amenable simple  ${\cal Z}$-stable \CA.
 Then 
 \beq
 K_0(A)_+=\{x\in K_0(A): \rho(A)>0\}\cup \{0\}\andeqn
 \Sigma(K_0(A))=\Sigma_\rho(K_0(A)).
 \eneq
In other words, 
$
(K_0(A), \Sigma(K_0(A)), \widetilde{T}(A),  \widetilde{\la e_A\ra}, \rho_A)
$ is
 a scaled simple ordered group pairing
(see Definition \ref{DElliott+} above and  Definitions 2.7 and 2.15 in \cite{GLrange})
(and $K_1(A)$ is a countable abelian group).
\end{thm}

The  range theorem  (the second ingredient)  below was first described  by
Elliott  in \cite{point-line} for finite case.

\begin{thm}[Gong-Lin, Theorem 5.3 of \cite{GLrange},  Theorem 5.6 of \cite{ERr}, and also \ref{range 0.35} above]
\label{GL-2}
Let
$(G_0, \Sigma_\rho(G_0), T, s, \rho)$
 be a simple ordered group pairing with distinguished set $\Sigma(G_0)$ and $G_1$  be a countable abelian group. Then there is a simple separable amenable ${\cal Z}$-stable \CA\, $A$ which satisfies the UCT such that
 \beq
 ((K_0(A), \Sigma_\rho(K_0(A)), \widetilde{T}(A), \widetilde{\la e_A\ra}, \rho_A), K_1(A))=((G_0, \Sigma_\rho(G_0), T, s, \rho), G_1).
 \eneq
 (i) In the case that $\widetilde{T}(A)\not=\{0\},$
A is unital if and only if $G_0$ has a unit $u,$ i.e.,  $u$ is the non-zero maximum in
$\Sigma_\rho(K_0(A))=\Sigma(K_0(A))$  and $\rho_A(u) = s$ (see Definition \ref{DElliott+}). If $\rho(G_0)\cap \Aff_+(T)\not=\{0\},$
then $K_0(A)_+=\{x\in K_0(A): \rho(x)>0\}\cup \{0\},$  $A$ can be chosen to have rational generalized tracial rank at most one (see Definition \ref{Dgtrf}) and be an inductive limit of subhomogeneous \CA s of spectra with dimension no more than $3.$ If
$\rho(G_0)\cap \Aff_+(T)=\{0\},$
 then A is stably projectionless and A can be chosen to
 be locally approximated by subhomogeneous \CA s with spectra having dimension no more than $3.$
If $G_1 = \{0\}$ and $G_0$ is torsion free, then $A$ can be chosen to be an inductive limit of 1-dimensional NCCW complexes with trivial $K_1.$

(ii) In the case that $\widetilde{T}(A)=\{0\},$ $A$ is purely infinite. $A$ is non-unital if and only if 
the distinguished set is empty, and unital if the distinguished set $\Sigma_\rho(K_0(A))=\{[1_A]\}.$

\end{thm}

Finally, let us present the isomorphism theorem:

\subsection{The isomorphism theorem}

We now present the third ingredient. 

\begin{thm}[Gong and Lin, Theorem 14.9 of \cite{GL-f}]\label{TisomorphismC}
Let $A$ and $B$ be finite separable  amenable simple  ${\cal Z}$-stable \CA s
satisfying the UCT.
Then $A\cong B$ if and only if
 \beq\label{TisomT-1.1}
 &&((K_0(A),  \Sigma(K_0(A)),  {\widetilde T}(A), {{\widehat{\la e_A\ra}}}, \rho_A), K_1(A))\\
 &&\hspace{0.3in}\cong ((K_0(B), \Sigma(K_0(B)), {\widetilde T}(B), {{\widehat{\la e_B\ra}}}, \rho_B), {{K_1(B)}}).
 \eneq
\end{thm}

(1)  If we restrict ourselves to  non-unital \CA s (i.e., we assume that
the \CA s are non-unital), the invariant used in Theorem \ref{TisomorphismC}
can be simplified to
$$((K_0(A),  {\widetilde T}(A), {{\widehat{\la e_A\ra}}}, \rho_A), K_1(A)).$$
This is because
\beq
\Sigma(K_0(A))=\{x\in K_0(A), x=[p]\,\, {\rm for\,\, some \,\,projection}\,\, p\in A\}\\
=\{x\in K_0(A): x\in K_0(A)_+\andeqn \rho_A(x)<{{\widehat{\la e_A\ra}}}\}
\eneq
which is determined by ${{\widehat{\la e_A\ra}}}.$

On the other hand, if we consider only the unital case, of course, the invariant may be
formulated   as $(K_0(A), [1_A], T(A), \rho_A, K_1(A)).$  However, we think
the invariant set should work for both unital and non-unital cases.
One then  should be able to tell from the invariant set itself
whether a simple \CA\, has a unit or not.
 The unified classification statements
 of Theorem \ref{TisomorphismC} and  Theorem \ref{GL-2} achieve that goal.
 
Note that Theorem \ref{GL-2} states that all Elliott invariants satisfying the description 
of Theorem \ref{GLrange1} can be achieved by a separable amenable simple ${\cal Z}$-stable 
\CA\, in  ${\cal N}$ (see Definition \ref{DUCT}). Omitting the UCT, but combining 
the purely infinite simple case (see \cite{Pclass} and \cite{KP}), let us state the following  
theorem:

\begin{thm}\label{Tfinal}
Let $A$ and $B$ be  separable  amenable simple  ${\cal Z}$-stable \CA s in ${\cal N}.$ 
Then $A\cong B$ if and only if
 \beq\label{TisomT-1.1}
 &&((K_0(A),  \Sigma_\rho(K_0(A)),  {\widetilde T}(A), {{\widehat{\la e_A\ra}}}, \rho_A), K_1(A))\\
 &&\hspace{0.3in}\cong ((K_0(B), \Sigma_\rho(K_0(B)), {\widetilde T}(B), {{\widehat{\la e_B\ra}}}, \rho_B), {{K_1(B)}}).
 \eneq
\end{thm}

Note that, in the purely infinite simple case, the Elliott invariants may be simplified as 
\beq
{\rm Ell}(A)=(K_0(A), [1_A], K_1(A)),\,\,\,{\rm or}\,\,\, {\rm Ell}(A)=(K_0(A), \emptyset,  K_1(A)).
\eneq

On the other hand, if $A$ is known to be stable, then 
the Elliott invariant may be simplified as 
\beq
{\rm Ell}(A)=(K_0(A),  \widetilde{T}(A), \rho_A, K_1(A))
\eneq
(finite or infinite). 


The Elliott program of classification of separable amenable \CA s does not end at Theorem \ref{GLrange1}, Theorem \ref{GL-2} and Theorem \ref{TisomorphismC}.
As mentioned earlier, there are known simple AH-algebras which are not ${\cal Z}$-stable.  There is also a question
whether non-simple ${\cal Z}$-stable \CA s could be classified, and if so, what would be the right formulation
of the invariant (see \cite{Ji-Jiang}, \cite{GJLP1}, \cite{GJLP2}, \cite{Jiang1}, \cite{CJL}, \cite{JLW}, \cite{GJL1}, \cite{GJL2}, and \cite{LNg}).
 Therefore  Theorem \ref{TisomorphismC} is not even the beginning of the end.
However, it does classify  the class of simple \CA s which  we have been mostly interested in for  decades.
The simple separable amenable \CA s classified by Theorem \ref{GLrange1}, Theorem \ref{GL-2}, and Theorem \ref{TisomorphismC} are of lower rank, or lower dimension
(rank zero or one).  A march to simple \CA s of higher rank has a long way to go (\cite{Vi} and \cite{Tex}). 
What is perhaps a modest beginning  was made recently in \cite{Elliott-Li-Niu}. Nevertheless, Theorem \ref{GLrange1}, Theorem \ref{GL-2} and Theorem \ref{TisomorphismC}
summarize  the decades of work and achievements  by many people in the Elliott program
for simple amenable ${\cal Z}$-stable \CA s.  It is, perhaps,  not an exaggeration  to state that
these three theorems (together with the classification of amenable purely infinite simple \CA s in \cite{Pclass}) represent the end of the beginning.

\vspace{0.2in}

Finally, it is exciting to note that a new approach
to the classification program---and in particular to
the finite case of Corollary 4.9---has recently been
developed (\cite{CGSTW}).


\vspace{0.2in}

{\bf Acknowledgements}: This paper is partially supported
by the Research Center for Operator Algebras at East China Normal University in Shanghai which is partially supported by Shanghai Key Laboratory of Pure Mathematics and Mathematical Practice, Science and Technology Commission of Shanghai Municipality (Grant No.~22DZ2229014).
The first  named author is supported by National Natural Science Foundation of China (Grant No.~11920101001), and the third named author is supported by a travel grant from the Simons Foundation (Grant No.~MP-TSM-00002606).

\noindent Guihua Gong\\
 University of Puerto Rico \\
San Juan, Puerto Rico \\
 USA, Rio Piedras, Puerto Rico 00936\\
ghgong@gmail.com

\vspace{0.2in}

\noindent Huaxin Lin \\
East China Normal University,\\
Shanghai, 200062\,\,\, and\\
University of Oregon \\
Eugene, Oregon 97402\\
USA\\
hlin@uoregon.edu

\vspace{0.2in}

\noindent Zhuang Niu\\
University of Wyoming \\
Laramie, Wyoming, 82072\\
USA\\
zniu@uwyo.edu

\end{document}